\documentclass{article}
\usepackage[utf8]{inputenc}
\usepackage{amsmath}
\usepackage{amsfonts}
\usepackage{amsthm}
\usepackage{enumitem}
\usepackage{hyperref}
\usepackage{tikz}
\usetikzlibrary{trees}
\usepackage{diagbox}
\usetikzlibrary{shapes.geometric}

\newtheorem{theorem}{Theorem}[section]
\newtheorem{definition}[theorem]{Definition}
\newtheorem{lemma}[theorem]{Lemma}
\newtheorem{proposition}[theorem]{Proposition}
\newtheorem{corollary}[theorem]{Corollary}
\newtheorem{conjecture}[theorem]{Conjecture}
\newtheorem{remark}[theorem]{Remark}
\newtheorem{notation}[theorem]{Notation}
\newtheorem{example}[theorem]{Example}

\usetikzlibrary{positioning, quotes}
\newdimen\nodeRad
\newdimen\nodeDist
\newdimen\nodeCondensedDist
\tikzset{
    position/.style args={#1:#2 from #3}{
        at=(#3.#1), anchor=#1+180, shift=(#1:#2)
    }
}

\tikzstyle{every node} = [circle, fill=gray!30,inner sep=0pt, minimum size = 11pt]
\tikzstyle{label} = [fill=white!0]
\tikzstyle{edge} = [draw=black]
\tikzstyle{edge_recent} = [draw=black, line width=2, very thick]
\tikzstyle{new} = [fill=red!30]
\tikzstyle{green} = [fill=green!30]
\tikzstyle{purple} = [fill=purple!30]
\tikzstyle{cyan} = [fill=cyan!30]

\usepackage{floatrow}
\newfloatcommand{capbtabbox}{table}[][\FBwidth]

\title{Catalan Number Sequences and Generalized Action Graphs}
\author{Drew Caldwell, Ali Cochran, Nathan Glisson, Bryce Jennings, \\
Katy McDicken, Luke Proctor, Sarah Klanderman, Amelia Tebbe}
\date{}

\begin{document}
\maketitle

\begin{abstract}
Action graphs emerged from work of Bergner and Hackney on category actions in the context of Reedy categories. Alvarez, Bergner, and Lopez showed that action graphs could be inductively generated without reference to category actions and have a close relationship with the sequence of Catalan numbers. These graphs were further generalized in work of Cressman, Lin, Nguyen, and Wiljanen, who showed that the Fuss-Catalan numbers have a similar relation to another set of inductively defined directed graphs. In our paper, we consider several other sequences related to the Catalan numbers, namely Catalan's triangle, $(a,b)$-Catalan numbers, internal triangles, and super Catalan numbers. We show action graphs cannot be generalized to Catalan's triangle, $(a,b)$-Catalan numbers, nor internal triangles. We also conjecture a method for constructing action graphs for the Super Catalan numbers. 
\end{abstract}

\section{Background}

 \textit{Action graphs} emerged from work of Bergner and Hackney on category actions in the context of Reedy categories \cite{Reedy}.  Alvarez, Bergner, and Lopez showed that action graphs could be inductively generated without reference to category actions, and they proved that the number of vertices added to $A_n$ is the $n$-th Catalan number \cite{ActionGraphs}.

\begin{definition}[\cite{ActionGraphs}]
The sequence of \textbf{action graphs} $\{A_n\}$ is defined inductively. Action graph $A_0$ is one vertex labeled 0 with no edges. Construct $A_{k+1}$ from $A_k$ by considering each vertex $v$ in $A_k$. For each path from $v$ to a vertex labeled $k$ in $A_k$, add a new edge with source $v$ to a new target vertex labeled $k+1$. Note that trivial paths are included in the paths considered from vertex $v$.
\end{definition}

\begin{figure}[h]
\begin{center}
\begin{tikzpicture}[> = stealth,shorten > = 1pt]
    \node[new] (root) at (0,0) {0};
\end{tikzpicture}
\hspace{1cm}
\begin{tikzpicture}[> = stealth,shorten > = 1pt]
    \node (root) at (0,0) {0};
    \node[new] (1) at (\nodeDist,0) {1};
    \path[->][edge_recent]  (root) to (1);
\end{tikzpicture}
\hspace{1cm}
\begin{tikzpicture}[> = stealth,shorten > = 1pt]
    \node (root) at (0,0) {0};
    \node (1) at (\nodeDist,0) {1};
    \node[new] (2) at (\nodeDist,\nodeDist) {2};
    \node[new] (3) at (0,\nodeDist) {2};
    \path[->][edge]  (root) to (1);
    \path[->][edge_recent]  (root) to (3);
    \path[->][edge_recent]  (1) to (2);
\end{tikzpicture}
\hspace{1cm}
\begin{tikzpicture}[> = stealth,shorten > = 1pt]
    \node (root) at (0,0) {0};
    \node (1) at (\nodeDist,0) {1};
    \node (2) at (\nodeDist,\nodeDist) {2};
    \node (3) at (0,\nodeDist) {2};
    
    \node[new] (4) at (0,2\nodeDist) {3};
    \node[new] (5) at (\nodeDist,2\nodeDist) {3};
    \node[new] (6) at (2\nodeDist,0) {3};
    \node[new] (7) at (-.7\nodeDist,.7\nodeDist) {3};
    \node[new] (8) at (-\nodeDist,0) {3};
    
    \path[->][edge]  (root) to (1);
    \path[->][edge]  (root) to (3);
    \path[->][edge]  (1) to (2);
    \path[->][edge_recent]  (3) to (4);
    \path[->][edge_recent]  (2) to (5);
    \path[->][edge_recent]  (1) to (6);
    \path[->][edge_recent]  (root) to (7);
    \path[->][edge_recent]  (root) to (8);
\end{tikzpicture}
\end{center}
\caption{Action Graphs $A_0$ through $A_3$, with new edges and vertices highlighted}
\label{actiongraphfigure}
\end{figure}
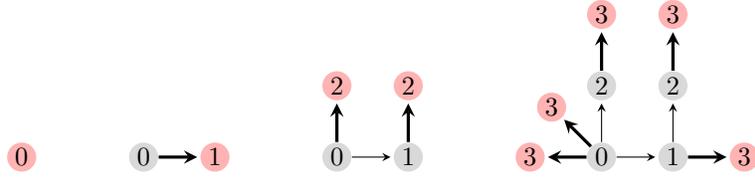

The first few graphs can be seen in Figure \ref{actiongraphfigure}.

\begin{definition}
    The \textbf{Catalan numbers} are a sequence of natural numbers given by 
$$ C_0 = 1, ~~C_{n} = \displaystyle\sum_{i=0}^{n-1} C_i C_{n-1-i} = \binom{2n}{n} \frac{1}{n+1}.$$
\end{definition}
The first few numbers of the sequence are $1,1,2,5,$ and $14$.

The action graphs in \cite{ActionGraphs} were further generalized in \cite{GenActionGraphs} by Cressman, Lin, Nguyen, and Wiljanen, who showed that the Fuss-Catalan numbers have a similar relation to another set of inductively defined directed graphs. The Fuss-Catalan numbers are a generalization of the Catalan numbers. 
    \begin{definition}[\cite{catalantext}, A14]
        The \textbf{Fuss-Catalan numbers} are defined by
    \[  C_{n,k} = \sum_{n_1+n_2+\cdots+n_{k+1}=n-1} \prod_{i=1}^{k+1} C_{n_i,k} = \frac{\binom{n(k+1)}{n}}{kn+1}.\]
    \end{definition}
   The sequence of the Fuss-Catalan numbers when $k=2$, that is $C_{n,2}$, is $1,1,3,12,55,\ldots$.
Observe that the Fuss-Catalan numbers agree with the Catalan numbers when $k=1$, that is $C_{n,1}=C_n$.

Cressman et al. expanded on the work of Alvarez, Bergner, and Lopez, developing new action graphs for the Fuss-Catalan numbers, which they called \textit{generalized action graphs} \cite{GenActionGraphs}. 

\begin{definition} \cite{GenActionGraphs} 
The \textbf{generalized action graphs} $T_{n,k}$ for the Fuss-Catalan numbers $C_{n,k}$ are defined inductively. For all $k$, $T_{0,k}$ is one vertex labeled $0$ with no edges. Construct $T_{n,k}$ from $T_{n-1,k}$ by considering each vertex $v$ in $T_{n-1,k}$. For each path of length $\ell$ from $v$ to a vertex labeled $n-1$ in $T_{n-1,k}$, add $\binom{\ell + k - 1}{\ell}$ new vertices labeled $n$ and new edges from $v$ to each of those new target vertices. 
\end{definition}

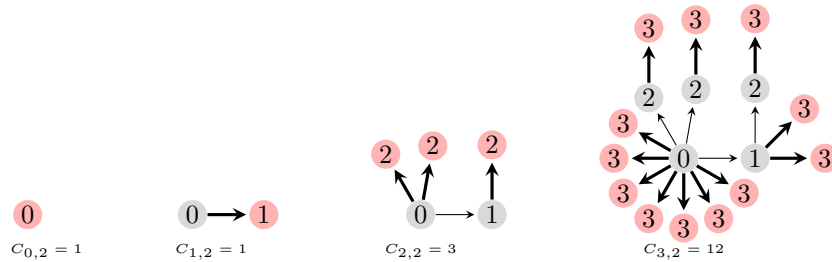
\begin{figure}[h]
\begin{center}

\begin{tikzpicture}[> = stealth,shorten > = 1pt]
    \node[new] (root) at (0,0) {0};
    \node[rectangle, fill=none] at (0.25, -.5) {\tiny{$C_{0,2}=1$}};
\end{tikzpicture}
\hspace{1cm}
\begin{tikzpicture}[> = stealth,shorten > = 1pt]
    \node (root) at (0,0) {0};
    \node[rectangle, fill=none] at (0.25, -.5) {\tiny{$C_{1,2}=1$}};
    \node[new] (1) at (\nodeDist,0) {1};
    \path[->][edge_recent]  (root) to (1);
\end{tikzpicture}
\hspace{1cm}
\begin{tikzpicture}[> = stealth,shorten > = 1pt]
    \node (root) at (0,0) {0};
    \node[rectangle, fill=none] at (0, -.5) {\tiny{$C_{2,2}=3$}};
    \node (1) at (\nodeDist,0) {1};
    \node[new] (2) [position=80:{\nodeRad} from root] {2};
    \node[new] (3) [position=120:{\nodeRad} from root] {2};
    \node[new] (4) [position=90:{\nodeRad} from 1] {2};
    \path[->][edge]  (root) to (1);
    \path[->][edge_recent]  (root) to (3);
    \path[->][edge_recent]  (root) to (2);
    \path[->][edge_recent]  (1) to (4);
\end{tikzpicture}
\hspace{1cm}
\begin{tikzpicture}[> = stealth,shorten > = 1pt]
    \node (root) at (0,0) {0};
    \node[rectangle, fill=none] at (0, -1.25) {\tiny{$C_{3,2}=12$}};
    \node (1) at (\nodeDist,0) {1};
    
    \node (2) [position=80:{\nodeRad} from root] {2};
    \node (3) [position=120:{\nodeRad} from root] {2};
    \node (4) [position=90:{\nodeRad} from 1] {2};
    
    \node[new] (5) [position=90:{\nodeRad} from 2] {3};
    \node[new] (6) [position=90:{\nodeRad} from 3] {3};
    \node[new] (7) [position=90:{\nodeRad} from 4] {3};
    \node[new] (8) [position=45:{\nodeRad} from 1] {3};
    \node[new] (9) [position=0:{\nodeRad} from 1] {3};
    
    \node[new] (16) [position=150:{\nodeRad} from root] {3};
    \node[new] (10) [position=180:{\nodeRad} from root] {3};
    \node[new] (11) [position=210:{\nodeRad} from root] {3};
    \node[new] (12) [position=240:{\nodeRad} from root] {3};
    \node[new] (13) [position=270:{\nodeRad} from root] {3};
    \node[new] (14) [position=300:{\nodeRad} from root]{3};
    \node[new] (15) [position=330:{\nodeRad} from root] {3};

    \path[->][edge]  (root) to (1);
    
    \path[->][edge]  (root) to (3);
    \path[->][edge]  (root) to (2);
    \path[->][edge]  (1) to (4);
    
    \path[->][edge_recent]  (2) to (5);
    \path[->][edge_recent]  (3) to (6);
    \path[->][edge_recent]  (1) to (8);
    \path[->][edge_recent]  (1) to (9);
    \path[->][edge_recent]  (4) to (7);
    
    \path[->][edge_recent]  (root) to (10);
    \path[->][edge_recent]  (root) to (11);
    \path[->][edge_recent]  (root) to (12);
    \path[->][edge_recent]  (root) to (13);
    \path[->][edge_recent]  (root) to (14);
    \path[->][edge_recent]  (root) to (15);
    \path[->][edge_recent]  (root) to (16);

\end{tikzpicture}
\end{center}
\caption{Generalized Action Graphs $T_{0,2}$ through $T_{3,2}$}
\end{figure}

Cressman et al. prove that the number of new vertices added to the generalized action graph $T_{n,k}$ is the Fuss-Catalan number $C_{n,k}$.

In this paper, we first discuss our notion of generalized action graphs. Then, in Sections \ref{Sec:AB}, \ref{Sec:Tri}, \ref{Sec:Internal}, and \ref{Sec:Super}, we consider several other sequences related to the Catalan numbers, namely $(a,b)-$Catalan numbers, Catalan's triangle, internal triangles, and super Catalan numbers respectively. We show action graphs cannot be generalized to Catalan's triangle, $(a,b)-$Catalan numbers, nor internal triangles. We also conjecture a method for constructing action graphs for the super Catalan numbers. \\

\subsection{Conventions and notation}

We use the following conventions for descendent and subtree in this paper. 

\begin{definition}
    A \textbf{descendent} of a vertex $v$ in a directed graph is any vertex $w$ such that there is a directed path from $v$ to $w$. For a rooted tree $T$ and a vertex $v$ of $T$, the \textbf{subtree of $T$ with root $v$} is the induced subgraph of vertex $v$ and its descendants. 
\end{definition}

Because generalized action graphs can easily become unwieldy in size, we introduce the following condensed notation to ease bookkeeping.

\begin{figure}[h]
\begin{center}
    \begin{tikzpicture}[> = stealth,shorten > = 1pt]
    \node (root) at (0,0) {0};
    \node[green,diamond] (1) at (\nodeDist,0) {1};
    \node[green,diamond] (2) at (-\nodeDist,0) {1};
    \path[->][edge]  (root) to (1);
    \path[->][edge]  (root) to (2);
    \node[new,rectangle] (3) at (1.7\nodeDist,.7\nodeDist) {2};
    \node[new,rectangle] (4) at (2\nodeDist,0) {2};
    \node[new,rectangle] (5) at (-2\nodeDist,0) {2};
    \node[new,rectangle] (6) at (-1.7\nodeDist, .7\nodeDist) {2};
    \path[->][edge_recent]  (1) to (3);
    \path[->][edge_recent]  (1) to (4);
    \path[->][edge_recent]  (2) to (5);
    \path[->][edge_recent]  (2) to (6);
    \node[new,rectangle] (7) at (.7\nodeDist, .7\nodeDist) {2};
    \node[new,rectangle] (8) at (-.7\nodeDist, .7\nodeDist) {2};
    \path[->][edge_recent]  (root) to (7);
    \path[->][edge_recent]  (root) to (8);
\end{tikzpicture}\hspace{.5in}
    \begin{tikzpicture}[> = stealth,shorten > = 1pt]
    \node (root) at (0,0) {0};
    \node[green,diamond] (1) at (\nodeCondensedDist,0) {1};
    \node[new,rectangle] (2) at (0,\nodeCondensedDist) {2};
    \node[new,rectangle] (4) at (\nodeCondensedDist,\nodeCondensedDist) {2};
{\footnotesize
    \draw[edge, ->] (root) -- (1) node[label, midway, above] {$\times 2$};
    \draw[edge, ->] (1) -- (4) node[label, midway, right] {$\times 2$};
    \draw[edge, ->] (root) -- (2) node[label, midway, right] {$\times 2$};
} 
    \path[->][edge] (root) to (1);
    \path[->][edge_recent] (root) to (2);
    \path[->][edge_recent] (1) to (4);
\end{tikzpicture}
\end{center}
\caption{A directed rooted tree in standard and condensed form}\label{fig:condensed}
\end{figure}
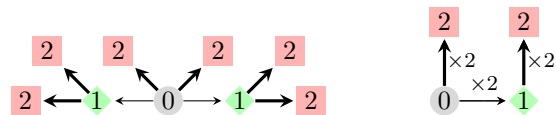

\begin{notation}\label{NotationCondensed}
     Since generalized action graphs have many identical subgraphs with the same labels, we collapse them. For an edge from a vertex labeled $k$ to a vertex labeled $n$, the multiplier, $\times m$, indicates the number of such edges from a vertex labeled $k$ to vertices labeled $n$ in the original graph. For a vertex in the condensed form, we can find the number of vertices it represents in the standard form by multiplying the labels along the path from the root to that vertex. For example, the upper right vertex labeled $2$ in the condensed graph in Figure \ref{fig:condensed} represents the $2\times 2=4$ vertices labeled $2$ that come off of the two vertices labeled 1 in the original graph.
\end{notation}

As in Figure \ref{fig:condensed}, we will sometimes change the color and shape of vertices in examples. This is intended only to highlight existing structure in the graphs, not to indicate any additional graph coloring or labeling. 

This paper will reference the sequences given in Table \ref{SequenceGlossary}.
\begin{table}[h]

\caption{
A table of sequences referenced in this paper.\label{SequenceGlossary}}
\begin{tabular}{|l|l|l|}
\hline
\textbf{Name}                   & \textbf{Notation \& Definition} & \textbf{Number Sequence} \\ \hline
Catalan numbers          & $C_{n} = \binom{2n}{n} \frac{1}{n+1}$                                     & $1,1,2,5,14,\ldots$                                                                                           \\ \hline
Fuss-Catalan numbers      & $C_{n,k} = \frac{1}{kn+1}\binom{n(k+1)}{n}$ & \begin{tabular}[c]{@{}l@{}} $k=1$: $1,1,2,5,14,\ldots$\\ $k=2$: $1,1,3,12,55,\ldots$ \end{tabular}     \\ \hline
Catalan's triangle        & $C(n,k) = \frac{n-k+1}{n+1} \binom{n+k}{k}$                                                                      & Table \ref{DrewTable}                                                                                     \\ \hline
$(a,b)-$Catalan numbers     & $Cat(a,b) = \frac{1}{a+b} \cdot \binom{a+b}{a}$                                                                                      & Table \ref{ABTable}                                                                                       \\ \hline
Weak $(a,b)-$Catalan numbers     & $cat(a,b) = \frac{1}{a+b} \cdot \binom{a+b}{a}$                                                                                      & Definition \ref{weak}                                                                                       \\ \hline
Internal triangles         & $t(n)=(n+2)C_{n-1}-2C_n$                                                                                                                                                              & $2,14,72, \ldots$                                                                                             \\ \hline
Super Catalan numbers      & $S(m,n) = \frac{(2m)!(2n)!}{m!n!(m+n)!}$                                                                                                                                                                        & \begin{tabular}[c]{@{}l@{}}$m=0$: $1, 2, 6, 20,70,$\ldots\\ $m=1$: $2,2,4,10,28,\ldots$\end{tabular} \\ \hline
\end{tabular}
\end{table}

\subsection{Generalized action graph properties}

We expand on the work of \cite{GenActionGraphs} by investigating whether similar graphs can be constructed for other sequences. Action graphs corresponding to the Catalan numbers and the Fuss-Catalan numbers have been defined, but there does not exist a general solution for any sequence of numbers in the literature. Based on the previous work on this topic, we define three axioms that describe generalized action graphs and explore some of the properties of such graphs.

\begin{definition}
    The sequence $\{G_n\}$ of \textbf{generalized action graphs} for a particular sequence $\{s_n\}$ of positive integers is a sequence of directed, labeled graphs such that:

\begin{enumerate}[label= Axiom \arabic{enumi}, leftmargin=*]
    \item\label{Axiom1} We define $G_0$ as the graph with $s_0$ vertices labeled 0 and no edges. We construct $G_n$ from $G_{n-1}$ by adding $s_n$ new vertices, which are each labeled $n$.
    \item\label{Axiom2} For vertex $v$ in $G_n$, the subtree of $G_n$ with root $v$ is isomorphic to some $G_k$ such that $k\leq n$. 
    \item\label{Axiom3} All leaves in the graph $G_n$ have label $n$.
\end{enumerate}
\end{definition}

Based on these axioms, some necessary conditions for generalized action graphs follow.
    \begin{lemma} \label{AG_One}
     In order for a sequence $\{s_n\}_{n\ge0}$ to form a valid sequence of generalized action graphs, it must have the property that $s_0 = 1$.   
    \end{lemma}
    \begin{proof} 
        Assume that $\{s_n\}$ has the property that $s_0\neq 1$. Then, $G_0$ would have multiple nodes labeled 0. If we consider the subtree rooted at any one of those 0 nodes in $G_0$, that subtree will contain just a single node, and therefore is not isomorphic to a generalized action graph in the sequence, which violates \ref{Axiom2}. Therefore, any sequence that has $s_0 \ne 1$ cannot have valid generalized action graphs.
    \end{proof}

 \begin{lemma}\label{squared}
    In order for a sequence $\{s_n\}_{n\geq0}$ to form a valid sequence of generalized action graphs, it must satisfy $s_2\geq s_{1}^2$.
    \end{lemma}
    \begin{proof}
        By Lemma \ref{AG_One}, $s_0 = 1$, and thus $G_0$ is a single node labeled 0. Suppose for contradiction that $s_1 = a$ and $s_2 = b$ such that $a^2 > b$. To create $G_1$, we add $s_1$ new leaves labeled 1. So $G_1$ has a total of $a$ new nodes labeled 1 with edges from the single node labeled 0. \ref{Axiom2} implies that in all future generalized action graphs, each node with path length zero to the newest node will result in $a$ new nodes. 
        
        To create $G_2$, we will add $s_2$ nodes labeled 2. However, as stated above, each node labeled 1 will result in $a$ new nodes labeled 2. This forces $a^2$ nodes labeled 2 to be added to the graph, which is greater than our allowed $s_2=b$ nodes labeled 2. Therefore, in order for a generalized action graph to exist, it must be that $s_2 \geq s_{1}^2$.
        \end{proof}

    The following example gives a more visual representation of the property described by Lemma \ref{squared}.

        \begin{example}
        Consider a sequence $\{a_n\}$ given by $1,4,14,48, \ldots$.

        Since $a_0 = 1$, this creates a single node on the first step labeled as 0. Then, since $a_1 = 4$, four nodes labeled with 1 have edges from the node labeled 0, which sets the precedent (by \ref{Axiom2}) that new nodes get four new leaves in the next step. Then, in the third step, since the most recently added nodes are known to receive four new nodes as descendants, we conclude that there should be at least 16 nodes labeled 2. However, $a_2 = 14$. If we add 16 nodes, then the graph will not satisfy \ref{Axiom1}.
    \end{example}
    \begin{center}
    \begin{tikzpicture}[> = stealth,shorten > = 1pt]
    \node[new] (root) at (0,0) {0};
\end{tikzpicture}
\hspace{.3cm}
\begin{tikzpicture}[> = stealth,shorten > = 1pt]
    \node (root) at (0,0) {0};
    \node[new] (1) [position=0:{\nodeRad} from root] {1};
    \node[new] (2) [position=60:{\nodeRad} from root] {1};
    \node[new] (3) [position=120:{\nodeRad} from root] {1};
    \node[new] (4) [position=180:{\nodeRad} from root] {1};
    \path[->][edge_recent]  (root) to (1);
    \path[->][edge_recent]  (root) to (2);
    \path[->][edge_recent]  (root) to (3);
    \path[->][edge_recent]  (root) to (4);
\end{tikzpicture}
\hspace{.3cm}
\begin{tikzpicture}[> = stealth,shorten > = 1pt]
    \node (root) at (0,0) {0};
    \node (1) [position=0:{\nodeRad} from root] {1};
    \node (2) [position=60:{\nodeRad} from root] {1};
    \node (3) [position=120:{\nodeRad} from root] {1};
    \node (4) [position=180:{\nodeRad} from root] {1};
    \node[new] (5) [position=0:{\nodeRad} from 1] {2};
    \node[new] (6) [position=20:{\nodeRad} from 1] {2};
    \node[new] (7) [position=40:{\nodeRad} from 1] {2};
    \node[new] (8) [position=60:{\nodeRad} from 1] {2};
    \node[new] (9) [position=40:{\nodeRad} from 2] {2};
    \node[new] (10) [position=60:{\nodeRad} from 2] {2};
    \node[new] (11) [position=80:{\nodeRad} from 2] {2};
    \node[new] (12) [position=100:{\nodeRad} from 2] {2};
    \node[new] (13) [position=80:{\nodeRad} from 3] {2};
    \node[new] (14) [position=100:{\nodeRad} from 3] {2};
    \node[new] (15) [position=120:{\nodeRad} from 3] {2};
    \node[new] (16) [position=140:{\nodeRad} from 3] {2};
    \node[new] (17) [position=120:{\nodeRad} from 4] {2};
    \node[new] (18) [position=140:{\nodeRad} from 4] {2};
    \node[new] (19) [position=160:{\nodeRad} from 4] {2};
    \node[new] (20) [position=180:{\nodeRad} from 4] {2};
    \path[->][edge]  (root) to (1);
    \path[->][edge]  (root) to (2);
    \path[->][edge]  (root) to (3);
    \path[->][edge]  (root) to (4);
    \path[->][edge_recent]  (1) to (5);
    \path[->][edge_recent]  (1) to (6);
    \path[->][edge_recent]  (1) to (7);
    \path[->][edge_recent]  (1) to (8);
    \path[->][edge_recent]  (2) to (9);
    \path[->][edge_recent]  (2) to (10);
    \path[->][edge_recent]  (2) to (11);
    \path[->][edge_recent]  (2) to (12);
    \path[->][edge_recent]  (3) to (13);
    \path[->][edge_recent]  (3) to (14);
    \path[->][edge_recent]  (3) to (15);
    \path[->][edge_recent]  (3) to (16);
    \path[->][edge_recent]  (4) to (17);
    \path[->][edge_recent]  (4) to (18);
    \path[->][edge_recent]  (4) to (19);
    \path[->][edge_recent]  (4) to (20);
\end{tikzpicture}
\end{center}

\section{$(a,b)-$Catalan numbers}\label{Sec:AB}

A sequence of numbers that is closely related the Catalan numbers is the sequence of $(a,b)-$Catalan numbers.
\begin{definition}\cite[A16]{catalantext}
For relatively prime integers $a \geq 0$ and $b \geq 0$, the \textbf{$(a,b)-$Catalan number}, denoted $Cat(a,b)$, is given by: \[ Cat(a,b) = \frac{1}{a+b} \cdot \binom{a+b}{a}=\frac{(a+b-1)!}{a!b!}. \] 
\end{definition}

Notice that the definition of an $(a,b)-$Catalan number requires that $a$ and $b$ be relatively prime integers. In order to assist in some of our exploration and calculations, it is useful to drop the relatively prime requirement. For non-negative integers $a$ and $b$, we define a generalization of the $(a,b)-$Catalan numbers that we call \emph{weak $(a,b)-$Catalan numbers}. As far as we are aware, this sequence has not previously been studied. 
\begin{definition} \label{weak}
    For any non-negative integers $a$ and $b$ that are not both zero, the \textbf{weak $(a,b)-$Catalan number} is given by \[ cat(a,b) = \frac{1}{a+b} \cdot \binom{a+b}{a}. \] 
\end{definition}
\begin{remark} \label{weakRemark}
 An important difference between standard $(a,b)-$Catalan numbers, $Cat(a,b)$, and weak $(a,b)-$Catalan numbers, $cat(a,b)$, is that a standard $(a,b)-$Catalan number will always be an integer, while a weak $(a,b)-$Catalan number may not be an integer. 
 \end{remark}
 
\begin{example}
Notice that $a=3$ and $b=6$ are not relatively prime numbers. We compute the weak $(a,b)-$Catalan number
\[ cat(3,6) = \frac{(3+6-1)!}{3!6!} = \frac{8!}{3!6!} = \frac{28}{3}. \]
\end{example}

For the remainder of this paper, when we refer to an $(a,b)-$Catalan number, we will mean the more general weak $(a,b)-$Catalan number.

\subsection{Properties of $(a,b)-$Catalan numbers}

We will discuss some important properties of $(a,b)-$Catalan numbers, including how they are connected to the Fuss-Catalan and Catalan numbers. 

Observe that for any $a$ and $b$, $cat(a,b) = cat(b,a)$ since 

\begin{equation}\label{ABProp1} 
cat(a,b) = \frac{(a+b-1)!}{a!b!} = \frac{(b+a-1)!}{b!a!} = cat(b,a).
\end{equation} 

Next, we will see the connection between the Fuss-Catalan numbers and the $(a,b)-$Catalan Numbers. Note that $cat(n,kn+1)=Cat(n,kn+1)$ here because $n$ and $kn+1$ are relatively prime.

\begin{theorem}\cite[A16a]{catalantext} \label{ABProp3}
 For all $n \geq 0$ and $k \in \mathbb{N}$,
 \[cat(n,kn+1) = C_{n,k}.\]
\end{theorem}

\begin{proof}

Recall that $C_{n,k} = \frac{1}{(kn+1)} \cdot \binom{n(k+1)}{n}$. Then, consider that
\[ C_{n,k}=\frac{1}{(kn+1)} \binom{n(k+1)}{n}=\frac{1}{kn+1}\cdot\frac{(n(k+1))!}{n!\cdot(n(k+1)-n)!}=
    \frac{(nk+n)!}{n!\cdot(nk+1)!}.\]
Now, notice that
$$cat(n,kn+1) = \frac{(n + (kn + 1) - 1)!}{n!(kn+1)!} = \frac{(nk + n)!}{n!(nk+1)!}.$$

Thus, $cat(n,kn+1) = C_{n,k}$ for all $n \geq 0$ and $k \in \mathbb{N}$.
\end{proof}

Corollary \ref{ABProp2} ties the $(a,b)-$Catalan numbers in with the Catalan numbers. This result can be observed using the previous theorem and the connection between the Fuss-Catalan numbers and Catalan numbers, $C_{n,1} = C_n$.

\begin{corollary}\cite[A16a]{catalantext}\label{ABProp2}
For all $n \geq 0$,
\[cat(n,n+1) = C_n.\]
\end{corollary}

To assist in proving Proposition \ref{diagonals}, we state the following lemma.

\begin{lemma} \label{ABProp4}
For $a \in \mathbb{N}$, $cat(a,1) = cat(1,a) = 1$.
\end{lemma}

\begin{proof}
Let $a \in \mathbb{N}$. By definition of $(a, b)$-Catalan number: \[ cat(a, 1) = \frac{(a+1 - 1)!}{a! \cdot 1!} = \frac{a!}{a!}=1. \]
Applying Equation \eqref{ABProp1}, $cat(a,1) = cat(1,a) = 1$.
\end{proof}

 Our goal is to explore whether we can make generalized action graphs for the $(a,b)-$Catalan numbers. The following proposition will assist us in that exploration.

\begin{proposition} \label{diagonals}
For all $n \geq 4$, $$n^2 \cdot cat(n+1, 1)^2 > n \cdot cat(n+2, 2).$$ 
\end{proposition}

\begin{proof}
We will prove this proposition using induction. 

For the base case, consider $n=4$. We have $4^2cat(5, 1)^2 > 4cat(6, 2)$ since $4^2cat
(5, 1)^2 = 16$ and $4cat(6, 2) = 14$.

Suppose there is some $k \geq 4$ such that $k^2 \cdot cat(k+1, 1)^2 > k \cdot cat(k+2, 2)$. 
By Lemma \ref{ABProp4}, we know that $cat(a,1) = 1$ for any $a \in \mathbb{N}$. So, we have $$k^2 > k \cdot cat(k+2, 2).$$ 

Since $k>0$, the above equation gives $k > cat(k+2, 2)$, which implies that $k \cdot \frac{k+4}{k+3} > cat(k+2, 2) \cdot \frac{k+4}{k+3}$. Note that by definition: \[ cat(k+2, 2) \cdot \frac{k+4}{k+3} = \frac{(k+3)!}{(k+2)!2!} \cdot \frac{k+4}{k+3} = \frac{(k+4)!}{(k+3)!2!} = cat(k+3, 2).\]
So, we have $k \cdot \frac{k+4}{k+3} > cat(k+3, 2)$.

Note that $m^2 + 4m + 3 > m^2 + 4m$ for any $m \in \mathbb{N}$. Since $$(m+1)(m+3) = m^2 + 4m + 3 > m^2 + 4m = m(m+4),$$ we can conclude that $m+1 > m \cdot \frac{m+4}{m+3}$. Thus the same is true for $k$, that is $$k+1 > k \cdot \frac{k+4}{k+3}>cat(k+3,2).$$

Since $cat(k+2, 1)^2 = 1^2 = 1$ by Lemma \ref{ABProp4}, we know $$k+1 = (k+1)\cdot cat(k+2,1)^2,$$ and therefore $(k+1) \cdot cat(k+2, 1)^2 > cat(k+3, 2)$, which implies $$(k+1)^2 \cdot cat(k+2, 1)^2 > (k+1) \cdot cat(k+3, 2).$$ This finishes the inductive step.

\end{proof}

The following two lemmas provide other properties of $(a,b)-$Catalan numbers.

\begin{lemma}
For any non-negative integers $a$ and $ b$, $$cat(a,b)=cat(a,b-1)+\frac{a-1}{a}cat(a-1,b).$$
\end{lemma}
\begin{proof}
We can use the definition of the $(a,b)-$Catalan numbers and simplify:

\begin{align*}
cat(a,b-1)+\frac{a-1}{a}cat(a-1,b)&=\frac{(a+b-2)!}{a!(b-1)!}+\frac{a-1}{a}\frac{(a+b-2)!}{(a-1)!b!}\\
&= \frac{b(a+b-2)!}{a!b!}+\frac{(a-1)(a+b-2)!}{a!b!}\\
&=\frac{(a+b-2)!(b+(a-1))}{a!b!}\\
&=\frac{(a+b-1)!}{a!b!}\\
&=\frac{1}{a+b}\frac{(a+b)!}{a!b!}\\
&=cat(a,b).
\end{align*}
\end{proof}

\begin{lemma}
For any integer $n\geq 2$, 
\[ cat(3,n)=cat(3,n-1)+\frac{n+1}{3}.\]
\end{lemma}
\begin{proof}
We will proceed by induction. As a base case, consider $n=2$. We have $$cat(3,2)=2=1+\frac{2+1}{3}=cat(3,2-1)+\frac{2+1}{3}.$$
Assume for some $k\geq 2$ that $cat(3,k)=cat(3,k-1)+\frac{k+1}{3}$.
Then consider, 

\begin{align*}
    cat(3,k+1)&=\frac{1}{3+k+1}\frac{(3+k+1)!}{3!(k+1)!} \tag{def of $(a,b)-$Catalan}\\
    &=\frac{k^2+5k+6}{3!}\\
    &=\frac{(k+1)k}{3!}+\frac{4k+6}{3!}\\
    &=\frac{(k+1)(k)}{3!}+\frac{2k+3}{3}\\
    &=\frac{1}{k+2}\frac{(k+2)!}{3!(k-1)!}+\frac{2k+3}{3}\\
    &=\frac{1}{3+k-1}\frac{(3+k-1)!}{3!(k-1)!}+\frac{2k+3}{3}\\
    &=cat(3,k-1)+\frac{k+1}{3}+\frac{k+2}{3}\tag{def of $(a,b)-$Catalan}\\
    &=cat(3,k)+\frac{k+2}{3}\tag{inductive hypothesis}\\
    &=cat(3,(k+1)-1) + \frac{(k+1)+1}{3}.
\end{align*}
Thus the inductive step holds, and our desired result follows.
\end{proof}

Observe that for arbitrary integers $a$ and $b$, the weak $(a,b)$-Catalan number $cat(a,b)$ is generally a fraction. However, action graphs are constructed based on whole number sequences. Therefore, we will introduce Catalan's triangle as a sequence that is closely related to the $(a,b)$-Catalan numbers, as described in Lemma \ref{TriToAB}, but whose entries are whole numbers.

\section{Catalan's triangle}\label{Sec:Tri}

  Another construction of numbers related to the Catalan numbers is Catalan's triangle.

\begin{definition}[\cite{OEISTriangle}]
    \textbf{Catalan's triangle} has entries denoted $C(n,k)$ that are defined as follows for integers $n$ and $k$  with $n\geq k\geq 0$:
    $$ C(n,k) = \frac{n-k+1}{n+1} \binom{n+k}{k}.$$
    
\end{definition}

We compute the entries of Catalan's triangle up to $n=8$ in Table \ref{DrewTable}.

\begin{table}[ht]
\caption{
A table of Catalan's triangle}
\label{DrewTable}
\begin{center}
\begin{tabular}{|l|l|l|l|l|l|l|l|l|l|}
\hline
\backslashbox{$n$}{$k$} & 0 & 1 & 2  & 3   & 4   & 5   & 6    & 7    & 8    \\ \hline
0   & 1 &   &    &     &     &     &      &      &      \\ \hline
1   & 1 & 1 &    &     &     &     &      &      &      \\ \hline
2   & 1 & 2 & 2  &     &     &     &      &      &      \\ \hline
3   & 1 & 3 & 5  & 5   &     &     &      &      &      \\ \hline
4   & 1 & 4 & 9  & 14  & 14  &     &      &      &      \\ \hline
5   & 1 & 5 & 14 & 28  & 42  & 42  &      &      &      \\ \hline
6   & 1 & 6 & 20 & 48  & 90  & 132 & 132  &      &      \\ \hline
7   & 1 & 7 & 27 & 75  & 165 & 297 & 429  & 429  &      \\ \hline
8   & 1 & 8 & 35 & 110 & 275 & 572 & 1001 & 1430 & 1430 \\ \hline
\end{tabular}
\end{center}
\end{table}

     We make the following observations about the entries of Catalan's triangle: 
     \begin{enumerate}
         \item $C(n,0) = 1$, for $n \ge 0.$
         \item $C(n,1) = n$, for $n \ge 1.$
         \item $C(n+1,k) = C(n+1,k-1)+C(n,k)$, for $1<k<n+1$.
         \item $C(n+1,n+1) = C(n+1,n)$, for $n \ge 1$.
     \end{enumerate}

     Because Catalan's triangle is an array rather than a sequence, there are many ways to create sequences from the entries. We consider building generalized action graphs for the sequences that arise as columns, rows, and diagonals of the triangle. We will describe these results in Section \ref{ab-AGs}, after we have shown the relationship between Catalan's triangle and the $(a,b)-$Catalan numbers.

\subsection{Relation to $(a,b)$-Catalan numbers}

        In this section, we will show that there is an interesting connection between the $(a,b)-$Catalan numbers and Catalan's triangle.

        \begin{lemma} \label{TriToAB}
            If $a \geq b\geq 0$ and $a\geq 1$, then $C(a-1,b)=(a-b)cat(a,b)$.
        \end{lemma}
        \begin{proof}
            Consider that for a Catalan's triangle entry $C(a-1,b)$, \begin{align*}
             C(a-1,b) &= \frac{a-b}{a} \binom{a+b-1}{b}\\ 
             &= \frac{a-b}{a} \cdot \frac{(a+b-1)!}{b!(a-1)!}\\ 
             &= \frac{(a-b)(a+b-1)!}{b!a!} \\
             &= (a-b)cat(a,b). \end{align*}
            
        \end{proof}

        This leads to some methods for concluding whether generalized action graphs can be formed for these sequences. First, we look at a consequence of Lemma \ref{TriToAB}. 
        \begin{lemma}
            Any $(a,b)-$Catalan number with $a\neq b$ can be written as a fraction with denominator $|a-b|$. In particular, when $a>b$,
            $$cat(a,b)=\frac{C(a-1,b)}{a-b}.$$
        \end{lemma}
        \begin{proof}
            Recall from Equation \eqref{ABProp1} that $cat(a,b)=cat(b,a)$. Without loss of generality, let $a>b$. By Lemma \ref{TriToAB}, $(a-b) \cdot cat(a,b) = C(a-1,b)$. Then $cat(a,b)=\frac{C(a-1,b)}{a-b}$, where $C(a-1,b)$ and $a-b$ are whole numbers. 
        \end{proof}
        To demonstrate this re-writing of $(a,b)-$Catalan numbers, we have Table \ref{ABTable}. Note that the diagonals represent where the differences in $a$ and $b$ are the same, and so they are written with the same denominators.

\begin{center}
\begin{table}[h]
\caption{
$(a,b)-$Catalan Numbers}
\label{ABTable}
\begin{tabular}{|l|l|l|l|l|l|l|l|l|l|}
\hline
\backslashbox{$a$}{$b$} & 0   & 1   & 2    & 3     & 4     & 5     & 6      & 7      & 8      \\
\hline
1   & 1/1 &     &      &       &       &       &        &        &        \\
\hline
2   & 1/2 & 1/1 &      &       &       &       &        &        &        \\
\hline
3   & 1/3 & 2/2 & 2/1  &       &       &       &        &        &        \\
\hline
4   & 1/4 & 3/3 & 5/2  & 5/1   &       &       &        &        &        \\
\hline
5   & 1/5 & 4/4 & 9/3  & 14/2  & 14/1  &       &        &        &        \\
\hline
6   & 1/6 & 5/5 & 14/4 & 28/3  & 42/2  & 42/1  &        &        &        \\
\hline
7   & 1/7 & 6/6 & 20/5 & 48/4  & 90/3  & 132/2 & 132/1  &        &        \\
\hline
8   & 1/8 & 7/7 & 27/6 & 75/5  & 165/4 & 297/3 & 429/2  & 429/1  &        \\
\hline
9   & 1/9 & 8/8 & 35/7 & 110/6 & 275/5 & 572/4 & 1001/3 & 1430/2 & 1430/1
\\ \hline
\end{tabular}
\end{table}
\end{center}

As previously noted, the weak $(a,b)-$Catalan numbers are frequently fractions rather than whole numbers. It is unclear to us what a fractional number of vertices should mean, so we will not try to directly construct generalized action graphs for this sequence. Instead, we will focus on the numerators in Table \ref{ABTable}, which, as a result of the preceding lemmas, are Catalan's triangle numbers. Observe that the numerators of Table \ref{ABTable} are exactly the entries in Table \ref{DrewTable}.

\subsection{Generalized action graphs for $(a,b)-$Catalan numbers and Catalan's triangle}\label{ab-AGs}

In this section we will find that generalized action graphs for the most part cannot be made for  sequences built from Catalan's triangle, and equivalently the numerators of the $(a,b)-$Catalan numbers. We will explore the columns, rows, and diagonals of Catalan's triangle. 

First, we will consider the columns.        
\begin{proposition}
    The columns of Catalan's triangle do not correspond to generalized action graphs, except in the case where $k=0$.
\end{proposition}
\begin{proof}
 The entries of the columns correspond to $C(n,k)$, where $k$ is fixed and $n$ is iterated. Except for the columns of $k=0$ and $k=1$, the columns have first entry greater than 1, since the $n=k$ diagonal is the Catalan numbers. (For illustration, see Table \ref{DrewTable}.) By Lemma \ref{AG_One}, these columns cannot be used to make a generalized action graph. For the column where $k=1$, notice that it contradicts Lemma \ref{squared}, and thus we cannot make a generalized action graph with this column. Therefore, the sequences made from the columns of Catalan's triangle cannot be modeled using (nontrivial) action graphs.
 
 For the column where $k=0$, we can make a somewhat trivial generalized action graph, where only one node is added each step.
 \end{proof}

Now, we will explore the rows. 
\begin{proposition}
    The rows of Catalan's triangle do not correspond to generalized action graphs.
\end{proposition}
\begin{proof}
    The entries of the rows correspond to $C(n,k)$, where $n$ is fixed and $k$ is iterated. First, we claim that for nonnegative $n$, $C(n,2) < C(n,1)^2$. To see this result, observe
    $$C(n,2) = \frac{n-2+1}{n+1} \binom{n+2}{2} = \frac{n-1}{n+1} \cdot \frac{(n+2)!}{2!((n+2)-2)!} =  \frac{n^2 + n - 2}{2}.$$
    Then by properties of Catalan's triangle, $C(n,1) = n$, thus $C(n,1)^2 = n^2$. Now, consider that
    \[ \frac{n^2 + n - 2}{2} - n^2 = \frac{-n^2 + n - 2}{2} \]
    and $\frac{-n^2 + n - 2}{2} < 0$ for all $n$. This implies that $$\frac{n^2 + n - 2}{2} - n^2 = C(n,2) - C(n,1)^2 < 0,$$ and thus $C(n,2) < C(n,1)^2$.
    
    So this means that using the rows to try and formulate generalized action graphs will lead to a contradiction of Lemma \ref{squared}, since the second value, $C(n,1)$, when squared is always more than the third value in the sequence, $C(n,2)$.
\end{proof}

Next, we will explore the diagonals. We will consider the $i$-th diagonal to be the sequence of entries $C(i+k,k)$ where $k$ is iterated. 

\begin{example}
    The third diagonal in Catalan's triangle is the sequence\\ $C(3,0), C(4,1), C(5,2), C(6,3), \ldots.$
\end{example}

\begin{proposition}
    Action graphs cannot be formulated for the $i$-th diagonal of Catalan's triangle when $i\geq 3$. 
\end{proposition}
\begin{proof}
Consider Proposition \ref{diagonals}, which says that $$m^2 \cdot cat(m+1, 1)^2 > m \cdot cat(m+2, 2)$$ for all $m \geq 4$. Combining this with Lemma \ref{TriToAB}, which implies that $$cat(m+j,j)=\frac{1}{m}C(m+j-1,j),$$ we see that 
\[m^2\frac{1}{m^2}C(m+1-1,1)^2>m\frac{1}{m}C(m+2-1,2),\]
and hence 
\[C((m-1)+1,1)^2>C((m-1)+2,2),\]
for $m-1\geq 3$. 
This results in a contradiction with Lemma \ref{squared}, which implies that we cannot construct generalized action graphs for the $i$-th diagonal of Catalan's triangle for any $i\geq 3$. 

\end{proof}

There still remain two Catalan's triangle diagonals to consider. The first diagonal consists of the Catalan numbers, and thus we know they have corresponding action graphs. Constructing generalized action graphs, or proving they do not exist, for the second diagonal of Catalan's triangle remains as future work.

\section{Internal triangles}\label{Sec:Internal}
One application of the Catalan numbers is counting the number of triangulations of polygons, using diagonals that do not intersect in the interior of the polygon. More specifically, for a convex $(n+2)$-gon $\mathcal{P}_{n+1}$, the number of triangulations is given by $C_n$ \cite{catalantext}. There is a sequence related to these triangulations referred to as the \textit{internal triangles}, which count the number of ways to draw triangles within a polygon without using the sides of the polygon.

\begin{lemma}[\cite{catalantext}, A5]
    The total number of internal triangles of a polygon with ${n+2}$ sides is given by $$t(n)=(n+2)C_{n-1}-2C_n = 2\binom{2n-3}{n-4}.$$
\end{lemma}

Via either of these formulas, we can compute $t(4) = 2, t(5) = 14, t(6)=72,$ and so on.

\begin{figure}[ht]
\begin{center}
\includegraphics[scale=0.250]{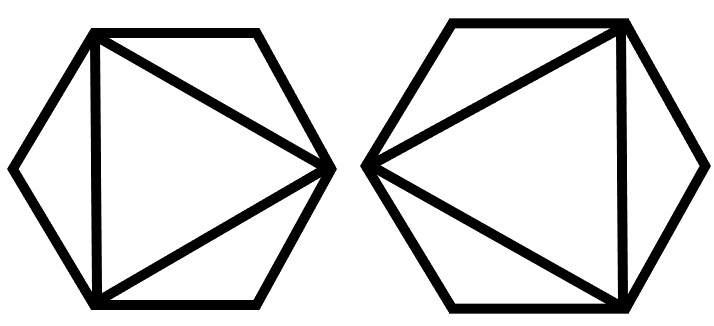}
\end{center}
\caption{$t(4)=2$, the number of internal triangles for hexagons.}
\label{internaltriangulation}
\end{figure}

We will use the term \textit{path length rule} to refer to the number of new vertices added to a vertex based on the length of paths from the vertex to a leaf labeled $n$. For $n\geq 0$, $\ell_n$ denotes the number of new vertices added for each path of length $n$.

\begin{theorem}
The internal triangles cannot be modeled using generalized action graphs. 
\end{theorem}
\begin{proof}
We will use path length rules to show that we cannot build consistent generalized action graphs beyond $t(6)$. Assume for a contradiction that we can model the internal triangle numbers using generalized action graphs. To satisfy Lemma \ref{AG_One}, we begin with one vertex labeled 0. We must construct the subsequent graphs to satisfy the three generalized action graph axioms.
 
Since $t(4)=2$, there is a single trivial path from the zero vertex to itself and we must have two new vertices labeled 1 coming off of our 0 vertex.  So, $\ell_0=2$.

In the next generalized action graph, our $\ell_0$ rule implies that we add two new vertices for each trivial path at a vertex labeled 1. Since $t(5)=14$, we then need to add ten additional vertices. Since there are two paths of length 1, we have $\ell_1=5$. These generalized action graphs are shown below in Figure \ref{internaltriangleactiongraph}.
 
\begin{figure}[ht]
\begin{center}
\begin{tikzpicture}[> = stealth,shorten > = 1pt]
    \node[new] (root) at (0,0) {0};
\end{tikzpicture}
\hspace{.75cm}
\begin{tikzpicture}[> = stealth,shorten > = 1pt]
    \node (root) at (0,0) {0};
    \node[new] (1) at (\nodeCondensedDist,0) {1};

    {\footnotesize
    \draw[edge_recent, ->] (root) -- (1) node[label, midway, below] {$\times 2$};
}
\end{tikzpicture}
\hspace{.75cm}
    \begin{tikzpicture}[> = stealth,shorten > = 1pt]
    \node (root) at (0,0) {0};
    \node (1) at (\nodeCondensedDist,0) {1};
    \node[new] (2) at (\nodeCondensedDist,\nodeCondensedDist) {2};
    \node[new] (3) at (0,\nodeCondensedDist) {2};
    {\footnotesize
    \draw[->,edge]  (root) -- (1) node[label, midway, below] {$\times 2$};
    \draw[->,edge_recent]  (root) -- (3) node[label, midway, right] {$\times 10$};
    \draw[->,edge_recent]  (1) -- (2) node[label, midway, right] {$\times 2$};
    }
    
\end{tikzpicture}
\hspace{.75cm}
\begin{tikzpicture}[> = stealth,shorten > = 1pt]
    \node (root) at (0,0) {0};
    \node (1) at (\nodeCondensedDist,0) {1};
    \node (2) at (\nodeCondensedDist,\nodeCondensedDist) {2};
    \node (3) at (0,\nodeCondensedDist) {2};
    
    \node[new] (4) at (0,2\nodeCondensedDist) {3};
    \node[new] (5) at (\nodeCondensedDist,2\nodeCondensedDist) {3};

    \node[new] (6) at (2\nodeCondensedDist,0) {3};

    \node[new] (8) at (-\nodeCondensedDist,0) {3};

    {\footnotesize
    \draw[->,edge]  (root) -- (1) node[label, midway, below] {$\times 2$};
    \draw[->,edge_recent]  (root) -- (3) node[label, midway, right] {$\times 10$};
    \draw[->,edge_recent]  (1) -- (2) node[label, midway, right] {$\times 2$};
    \draw[->,edge_recent]  (3) -- (4) node[label, midway, right] {$\times 2$};
    \draw[->,edge_recent]  (2) -- (5) node[label, midway, right] {$\times 2$};
    \draw[->,edge_recent]  (root) -- (8) node[label, midway, below] {$\times 50$};
    \draw[->,edge_recent]  (1) -- (6) node[label, midway, below] {$\times 10$};
    }
\end{tikzpicture}
\end{center}
\caption{Attempted generalized action graphs for the internal triangle numbers}
\label{internaltriangleactiongraph}
\end{figure}
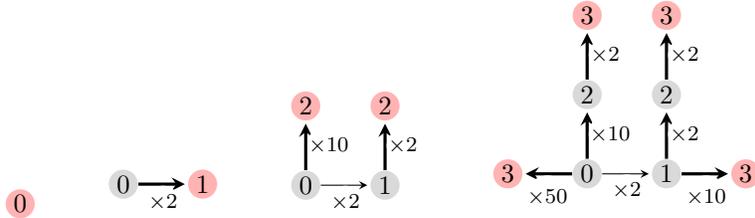

    To create the next generalized action graph, note that we have fourteen vertices labeled 2 with trivial edges. By the path length rule $\ell_0$, we add two vertices labeled 3 for each vertex labeled 2, giving us 28 new vertices total. We have four paths of length one starting at vertices labeled 1, and ten paths of length one starting at vertex 0. Since $\ell_1 = 5$, we add $4\times 5 = 20$ new vertices with edges from the vertices labeled 1 and $10\times 5 = 50$ new vertices with edges from vertex 0. Recall that $t(6)=72$, but we already added ninety-eight ($28+20+50$) new vertices to this generalized action graph. Because of this contradiction, the internal triangle numbers cannot be modeled using generalized action graphs.
     
\end{proof}

\section{The super Catalan numbers}\label{Sec:Super}
The super Catalan numbers are a generalized form of the Catalan numbers with two arguments, $m$ and $n$.  

\begin{definition}[\cite{catalantext}, A17]\label{SuperCatalnDef}
    The \textbf{super Catalan numbers} are defined by $$S(m,n) = \frac{(2m)!(2n)!}{m!n!(m+n)!}.$$
\end{definition}

When $m = 0$, the first few numbers in the sequence are $1, 2, 6, 20,$ and $ 70$.

We wish to construct generalized action graphs for the sequence of super Catalan numbers where $m=0$. In the rest of this section, we will state a method for constructing a sequence of directed graphs, conjecture that these graphs are indeed generalized action graphs for the sequence $S(0,n)$, and discuss our progress towards proving this conjecture. 

\begin{definition}\label{SuperCatAGConstruction}
    We construct the sequence \textbf{generalized action graphs}, denoted $\{G_n\}$, for the super Catalan numbers as sequence of directed graphs defined inductively in the following way. The graph $G_0$ is a single vertex labeled 0. To construct $G_{n+1}$ from $G_n$, consider each vertex $v$ in $G_n$. For each $0\leq\ell\leq n$, add $\displaystyle p(v,\ell)\cdot\frac{2}{2^\ell}$ new vertices labeled $n+1$ with edges from $v$, where $p(v,\ell)$ is the number of paths of length $\ell$ from $v$ to vertices labeled $n$ in $G_n$. 
\end{definition}

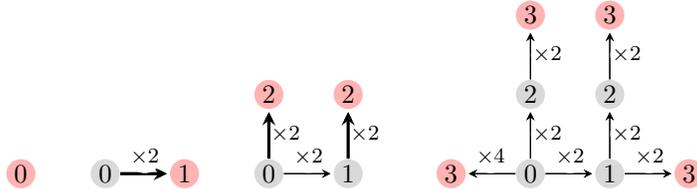
\begin{figure}[ht]
\begin{center}
\begin{tikzpicture}[> = stealth,shorten > = 1pt]
    \node[new] (root) at (0,0) {0};
\end{tikzpicture}
\hspace{0.5cm}
\begin{tikzpicture}[> = stealth,shorten > = 1pt]
     \node (root) at (0,0) {0};
    \node[new] (1) at (\nodeCondensedDist,0) {1};
    \path[->][edge_recent]  (root) to (1);

    {\footnotesize
    \draw[edge, ->] (root) -- (1) node[label, midway, above] {$\times 2$};
    }
\end{tikzpicture}
\hspace{0.5cm}
\begin{tikzpicture}[> = stealth,shorten > = 1pt]
    \node (root) at (0,0) {0};
    \node (1) at (\nodeCondensedDist,0) {1};
    \node[new] (2) at (0,\nodeCondensedDist) {2};
    \node[new] (4) at (\nodeCondensedDist,\nodeCondensedDist) {2};
{\footnotesize
    \draw[edge, ->] (root) -- (1) node[label, midway, above] {$\times 2$};
    \draw[edge, ->] (1) -- (4) node[label, midway, right] {$\times 2$};
    \draw[edge, ->] (root) -- (2) node[label, midway, right] {$\times 2$};
} 
    \path[->][edge] (root) to (1);
    \path[->][edge_recent] (root) to (2);
    \path[->][edge_recent] (1) to (4);
\end{tikzpicture}
\hspace{0.5cm}
\begin{tikzpicture}[> = stealth,shorten > = 1pt]
    \node (root) at (0,0) {0};
    \node (1) at (\nodeCondensedDist,0) {1};
    \node (2) at (0,\nodeCondensedDist) {2};
    \node[new] (3) at (0,2\nodeCondensedDist) {3};
    \node (4) at (\nodeCondensedDist,\nodeCondensedDist) {2};
    \node[new] (5) at (\nodeCondensedDist,2\nodeCondensedDist) {3};
    \node[new] (6) at (-\nodeCondensedDist,0) {3};
    \node[new] (7) at (2\nodeCondensedDist,0) {3};

    \path[->][edge] (root) to (1);
    \path[->][edge] (root) to (2);
    \path[->][edge] (2) to (3);
    \path[->][edge] (1) to (4);
    \path[->][edge] (4) to (5);
    \path[->][edge] (root) to (6);
    \path[->][edge] (1) to (7);

{\footnotesize
    \draw[edge, ->] (root) -- (1) node[label, midway, above] {$\times 2$};
    \draw[edge, ->] (1) -- (4) node[label, midway, right] {$\times 2$};
    \draw[edge, ->] (4) -- (5) node[label, midway, right] {$\times 2$};
    \draw[edge, ->] (root) -- (2) node[label, midway, right] {$\times 2$};
    \draw[edge, ->] (2) -- (3) node[label, midway, right] {$\times 2$};
    \draw[edge, ->] (root) -- (6) node[label, midway, above] {$\times 4$};
    \draw[edge, ->] (1) -- (7) node[label, midway, above] {$\times 2$};
}
\end{tikzpicture}
\end{center}
\caption{Generalized action graphs $G_0$ through $G_3$ for super Catalan numbers}
\label{SuperCatalanGraphs}
\end{figure}

For an example of the construction of $G_0$, $G_1$, $G_2$, and $G_3$ see Figure \ref{SuperCatalanGraphs}. Note that we use the condensed notation introduced in Notation \ref{NotationCondensed}.

Definition \ref{SuperCatAGConstruction} satisfies \ref{Axiom3}, since each leaf in $G_n$ will be labeled $n$ and will have a path of length 0 to itself, and thus will have at least two new descendants, which will be labeled $n+1$, in $G_{n+1}$.

\begin{conjecture}\label{SuperCatAGConjecture}
     Each graph $G_n$ in the sequence of graphs constructed in Definition \ref{SuperCatAGConstruction} satisfies \ref{Axiom1} and \ref{Axiom2}.
\end{conjecture}

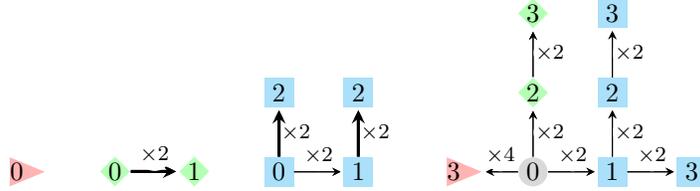
\begin{figure}[ht]
\begin{center}
\begin{tikzpicture}[> = stealth,shorten > = 1pt]
    \node[new, isosceles triangle] (root) at (0,0) {0};
\end{tikzpicture}
\hspace{0.5cm}
\begin{tikzpicture}[> = stealth,shorten > = 1pt]
     \node[green,diamond] (root) at (0,0) {0};
    \node[green,diamond] (1) at (\nodeCondensedDist,0) {1};
    \path[->][edge_recent]  (root) to (1);

    {\footnotesize
    \draw[edge, ->] (root) -- (1) node[label, midway, above] {$\times 2$};
    }
\end{tikzpicture}
\hspace{0.5cm}
\begin{tikzpicture}[> = stealth,shorten > = 1pt]
    \node[cyan,rectangle] (root) at (0,0) {0};
    \node[cyan,rectangle] (1) at (\nodeCondensedDist,0) {1};
    \node[cyan,rectangle] (2) at (0,\nodeCondensedDist) {2};
    \node[cyan,rectangle] (4) at (\nodeCondensedDist,\nodeCondensedDist) {2};
{\footnotesize
    \draw[edge, ->] (root) -- (1) node[label, midway, above] {$\times 2$};
    \draw[edge, ->] (1) -- (4) node[label, midway, right] {$\times 2$};
    \draw[edge, ->] (root) -- (2) node[label, midway, right] {$\times 2$};
} 
    \path[->][edge] (root) to (1);
    \path[->][edge_recent] (root) to (2);
    \path[->][edge_recent] (1) to (4);
\end{tikzpicture}
\hspace{0.5cm}
\begin{tikzpicture}[> = stealth,shorten > = 1pt]
    \node (root) at (0,0) {0};
    \node[cyan,rectangle] (1) at (\nodeCondensedDist,0) {1};
    \node[green,diamond] (2) at (0,\nodeCondensedDist) {2};
    \node[green,diamond] (3) at (0,2\nodeCondensedDist) {3};
    \node[cyan,rectangle] (4) at (\nodeCondensedDist,\nodeCondensedDist) {2};
    \node[cyan,rectangle] (5) at (\nodeCondensedDist,2\nodeCondensedDist) {3};
    \node[new,isosceles triangle] (6) at (-\nodeCondensedDist,0) {3};
    \node[cyan,rectangle] (7) at (2\nodeCondensedDist,0) {3};

    \path[->][edge] (root) to (1);
    \path[->][edge] (root) to (2);
    \path[->][edge] (2) to (3);
    \path[->][edge] (1) to (4);
    \path[->][edge] (4) to (5);
    \path[->][edge] (root) to (6);
    \path[->][edge] (1) to (7);

{\footnotesize
    \draw[edge, ->] (root) -- (1) node[label, midway, above] {$\times 2$};
    \draw[edge, ->] (1) -- (4) node[label, midway, right] {$\times 2$};
    \draw[edge, ->] (4) -- (5) node[label, midway, right] {$\times 2$};
    \draw[edge, ->] (root) -- (2) node[label, midway, right] {$\times 2$};
    \draw[edge, ->] (2) -- (3) node[label, midway, right] {$\times 2$};
    \draw[edge, ->] (root) -- (6) node[label, midway, above] {$\times 4$};
    \draw[edge, ->] (1) -- (7) node[label, midway, above] {$\times 2$};
}
\end{tikzpicture}
\end{center}
\caption{Highlighting the isomorphic subtrees of $G_3$}
\label{subgraphfigure}
\end{figure}

In support of this conjecture, consider the following example. In Figure \ref{subgraphfigure}, we use colors to show that each subtree of $G_3$ is isomorphic to a previous generalized action graph. Because the graphs are drawn in condensed notation, the multiplier on the edge from the $0$ vertex indicates the number of subtrees that are isomorphic to a previous generalized action graphs. In particular, considering the subtrees rooted at vertices adjacent to the zero vertex: $G_3$ has two subtrees isomorphic to $G_2$, two subtrees isomorphic to $G_1$, and four subtrees isomorphic to $G_0$, each with the labels of their vertices shifted up by an appropriate amount.

We will discuss our progress towards proving Conjecture \ref{SuperCatAGConjecture}.
The data found in $\{G_n\}$ can be encoded into tables, which we call \textit{$n$-tables}. 

\begin{definition}\label{ntable}
    Let $K_{\ell,v,n}$ be the number of paths of length $\ell$ in $G_n$ that start at a vertex labeled $v$ and end at a vertex labeled $n$. For a given $n$, the table of $K_{\ell,v,n}$ for all values of $\ell$ and $v$ is called the \textbf{$n$-table}.
    
\end{definition}

   Figure \ref{superactiongraph} shows $G_3$, the conjectured generalized action graph for $S(0,3)$ in condensed notation, and Table \ref{pathsaddedtable} is the corresponding $n$-table. For example, to find $K_{1,2,3}$, count the paths of length 1 from vertices labeled 2 to vertices labeled 3 by taking the $2 \times$ \textcolor{blue}{$2=4$} paths from the leftmost 2 to vertices labeled 3, and adding the $2 \times 2 \times$ \textcolor{blue}{$2=8$} paths from the rightmost 2 to vertices labeled 3. Together that is \textcolor{blue}{12}, as seen in the corresponding entry of Table \ref{pathsaddedtable}.

\begin{figure}[ht]
\begin{floatrow}
\ffigbox[.42\textwidth]{%
  \begin{tikzpicture}[> = stealth,shorten > = 1pt]
    \node(root) at (0,0) {0};
    \node (1) at (\nodeCondensedDist,0) {1};
    \node (2) at (0,\nodeCondensedDist) {2};
    \node (3) at (0,2\nodeCondensedDist) {3};
    \node (4) at (\nodeCondensedDist,\nodeCondensedDist) {2};
    \node (5) at (\nodeCondensedDist,2\nodeCondensedDist) {3};
    \node (6) at (-\nodeCondensedDist,0) {3};
    \node (7) at (2\nodeCondensedDist,0) {3};
{\footnotesize
    \draw[edge, ->] (root) -- (1) node[label, midway, above] {$\times 2$};
    \draw[edge, ->] (1) -- (4) node[label, midway, right] {$\times 2$};
    \draw[edge, ->] (4) -- (5) node[label, midway, right] {\textcolor{blue}{$\times 2$}};
    \draw[edge, ->] (root) -- (2) node[label, midway, right] {$\times 2$};
    \draw[edge, ->] (2) -- (3) node[label, midway, right] {\textcolor{blue}{$\times 2$}};
    \draw[edge, ->] (root) -- (6) node[label, midway, above] {$\times 4$};
    \draw[edge, ->] (1) -- (7) node[label, midway, above] {$\times 2$};
}
    \path[->][edge] (root) to (1);
    \path[->][edge] (root) to (2);
    \path[->][edge] (2) to (3);
    \path[->][edge] (1) to (4);
    \path[->][edge] (4) to (5);
    \path[->][edge] (root) to (6);
    \path[->][edge] (1) to (7);
    \end{tikzpicture}}%
{%
  \caption{Conjectured $G_3$}\label{superactiongraph}%
}
\capbtabbox{%
  \begin{tabular}{|l|l|l|l|l|}
    \hline
    \textbf{$n$ = 3} & \textbf{$v$ = 0} & \textbf{$v$ = 1} & \textbf{$v$ = 2} & \textbf{$v$ = 3} \\ \hline
    \textbf{$\ell$ = 0} & 0     & 0     & 0     & 20     \\ \hline
    \textbf{$\ell$ = 1} & 4     & 4     & \textcolor{blue}{12}    & 0      \\ \hline
    \textbf{$\ell$ = 2} & 8     & 8     & 0     & 0      \\ \hline
    \textbf{$\ell$ = 3} & 8     & 0     & 0     & 0      \\ \hline
    
    \end{tabular}
        
}{%
   \caption{$K_{\ell, v, 3}$ values for $G_3$}%
   \label{pathsaddedtable}
}
\end{floatrow}
\end{figure}

Each entry in the $n$-table is completely determined by the graph $G_n$ that it represents. Assuming Conjecture \ref{klemma} is true, we found a way to compute the entries without having to count each path from each physically drawn generalized action graph:
\begin{itemize}
    \item When $\ell = 0$, $K_{0, n, v} = 0$ for all $v$ values except when $v = n$; in that case $K_{0, n, n} = S(0,n)$. 
    \item For $\ell>0$, entry $K_{l,v,n}$ in the $n$-table can be found using Conjecture \ref{klemma}, which uses the entries from the corresponding column of the $(n-1)$-table for the graph $G_{n-1}$.
\end{itemize}

\begin{conjecture}\label{klemma}
    The number of paths of length $\ell$ in $G_{n+1}$ that start at a vertex labeled $v$ and end at a vertex labeled $n+1$ is given by
    $$K_{\ell,v,n+1} = \displaystyle \sum_{i=0}^{n+1-\ell-v} \frac{2}{2^i} K_{\ell-1+i,v,n}.$$
\end{conjecture}

We have checked that Conjecture \ref{klemma} holds up through the $7$-tables, but for larger $n$ the generalized action graphs and the associated tables become cumbersome to both create and count. As an example of how to use Conjecture \ref{klemma}, consider the $n$-tables for the graphs of $S(0,3)$ and $S(0,4)$ in Tables \ref{pathsaddedtable3} and \ref{pathsaddedtable4} respectively.

\begin{minipage}{.45\textwidth}
    \bigskip
    
    \begin{table}[H]

    \begin{tabular}{|l||l|l|l|l|}
    \hline
    \backslashbox{$\ell$}{$v$} &  0 &  1 &  2 &  3 \\ 
\hline\hline
    $\ell$ = 0 & \textcolor{violet}{0}     & 0     & 0     & 20     \\ \hline
    $\ell$ = 1 & \textcolor{violet}{4}     & 4     & 12    & 0      \\ \hline
     $\ell$ = 2 & \textcolor{violet}{8}     & 8     & 0     & 0      \\ \hline
    $\ell$ = 3 & \textcolor{violet}{8}     & 0     & 0     & 0      \\ \hline
    \end{tabular}
         \caption{$K_{\ell, v, 3}$ values for $G_3$}
         \label{pathsaddedtable3}
    \end{table}
    \end{minipage}
    \begin{minipage}{.45\textwidth}
    \hspace{1cm}
    
        \begin{table}[H]
    
    \begin{tabular}{|l||l|l|l|l|l|}
    \hline
    \backslashbox{$\ell$}{$v$}&  0 &  1 & 2 &  3 &  4 \\ \hline \hline
    $\ell$ = 0 & 0     & 0     & 0     & 0     & 70     \\ \hline
    $\ell$ = 1 & \textcolor{violet}{10}     & 8     & 12    & 40   & 0      \\ \hline
    $\ell$ = 2 & 20     & 16     & 24     & 0      & 0      \\ \hline
    $\ell$ = 3 & 24     & 16     & 0     & 0       & 0     \\ \hline
    $\ell$ = 4 & 16    & 0     & 0     & 0     & 0 \\ \hline 
    \end{tabular}
         \caption{$K_{\ell, v, 4}$ values for $G_4$}
         \label{pathsaddedtable4}
    \end{table}
    \end{minipage}
    \bigskip

We will show how to compute $K_{1, 0, 4}= \textcolor{violet}{10}$ from Table \ref{pathsaddedtable4} using Conjecture \ref{klemma}:

\begin{align*}
    K_{1,0,4} &= \sum_{i=0}^{3} \frac{2}{2^i} K_{i,0,3}\\
    &= \frac{2}{2^0}(\textcolor{violet}{0})+ \frac{2}{2^1}(\textcolor{violet}{4})+ \frac{2}{2^2}(\textcolor{violet}{8})+\frac{2}{2^3}(\textcolor{violet}{8})\\
    &= \textcolor{violet}{10}.
\end{align*}

If Conjecture \ref{klemma} holds, we believe that Conjecture \ref{SKconjecture} can be used to compute the next super Catalan number, $S(0,n+1)$, using the $K_{\ell,v,n}$ values from the $n$-table for the graph of $S(0, n)$, which would show the graph $G_n$ satisfies \ref{Axiom1}.

\begin{conjecture}\label{SKconjecture}
The subsequent super Catalan number can be computed from the $n$-table of its previous action graph via
    $$S{(0,n+1)} = \displaystyle \sum_{\ell=0}^{n} \left(\frac{2}{2^\ell} \sum_{v=0}^{n} K_{\ell,v,n}\right).$$
\end{conjecture}

 This conjecture uses the sums of each column from Table \ref{pathsaddedtable3}. We will check that Conjecture \ref{SKconjecture} is true for $n=3$:

\begin{align*}
    S(0,4) &= \sum_{\ell=0}^{3} \left(\frac{2}{2^\ell} \sum_{v=0}^{3} K_{\ell,v,3}\right)\\
    &= \frac{2}{2^0}\sum_{v=0}^{3}K_{0,v,3} + \frac{2}{2^1}\sum_{v=0}^{3}K_{1,v,3} + \frac{2}{2^2}\sum_{v=0}^{3}K_{2,v,3} + \frac{2}{2^3}\sum_{v=0}^3K_{3,v,3}\\
    &= 2(0+0+0+20) + 1(4+4+12+0) + \frac{1}{2}(8+8+0+0) + \frac{1}{4}(8+0+0+0)\\
     &= 70,
\end{align*}
which is known to be the super Catalan number $S(0, 4)$.

\section{Acknowledgements}
This project was supported by a grant from the Center for Undergraduate Research in Mathematics, which is funded by the National Science Foundation DMS awards 0636648, 1148695, and 1722563. Thanks also to Finn Richert and Maddie Dobrocky for their work in the early stages of these projects.


\begin{thebibliography}{99} 

\bibitem{catalantext} Richard Stanley. ``Catalan Numbers'' \textit{Cambridge University Press, New York} (2015)


\bibitem{ActionGraphs} Alvarez, Bergner, Lopez. ``Action Graphs and Catalan Numbers.'' \textit{Journal of Integer Sequences} Vol. 18 (2015)
\href{https://arxiv.org/abs/1503.00044v1}{https://arxiv.org/abs/1503.00044v1}

\bibitem{GenActionGraphs} Danielle Cressman, Jonathan Lin, An Nguyen, and Luke Wiljanen. 
``Generalized action graphs.'' (In preparation) 

\bibitem{Reedy}
Julia E. Bergner and Philip Hackney. ``Reedy categories which encode the notion of category actions.''
\textit{Fundamenta Mathematicae} 228.3 (2015) p. 193-222.
\href{http://eudml.org/doc/282637}{http://eudml.org/doc/282637}

 \bibitem{Gould}
 H.W.Gould and Jocelyn Quaintance.
 ``Combinatorial Identities: Table II: Advanced Techniques for Summing Finite Series.'' \href{https://math.wvu.edu/~hgould/Vol.5.PDF}{https://math.wvu.edu/~hgould/Vol.5.PDF}

\bibitem{OEISTriangle} Online Encyclopedia of Integer Sequences,
\href{https://oeis.org/A009766}{https://oeis.org/A009766}
\end{thebibliography}
\end{document}